\documentclass[12pt]{amsart}
\usepackage{amssymb}
\newtheorem{theorem}{Theorem}[section]
\newtheorem{lemma}[theorem]{Lemma}
\newtheorem{prop}[theorem]{Proposition}
\newtheorem{cor}[theorem]{Corollary}

\newtheorem{notn}[theorem]{Notation}
\theoremstyle{definition}
\newtheorem{example}[theorem]{Example}
\newtheorem{remark}[theorem]{Remark}
\newtheorem{definition}[theorem]{Definition}

\newtheorem{question}[theorem]{Question}

\newcommand{\R}{\mathbb{R}}
\newcommand{\Z}{\mathbb{Z}}
\newcommand{\N}{\mathbb{N}}
\newcommand{\Q}{\mathbb{Q}}

\newcommand{\e}{\mathbf{e}}
\newcommand{\lam}{\ensuremath{{\lambda}}}
\newcommand{\NP}{\mathrm{NP}}

\newcommand{\conv}{\mathrm{conv}}

\newcommand{\xvec}[1]{\ensuremath{x_{1}, \ldots, x_{#1}}}
\newcommand{\lcm}{\ensuremath{\mathrm{lcm}}}
\newcommand{\grp}{\ensuremath{\mathrm{grp}}}

\newcommand{\blam}{\boldsymbol{\lambda}}
\newcommand{\blambda}{\boldsymbol{\lambda}}
\newcommand{\bomega}{\boldsymbol{\omega}}
\newcommand{\bbeta}{\boldsymbol{\beta}}
\newcommand{\balpha}{\boldsymbol{\alpha}}
\newcommand{\bgamma}{\boldsymbol{\gamma}}
\newcommand{\bdelta}{\boldsymbol{\delta}}
\newcommand{\bnu}{\boldsymbol{\nu}}
\newcommand{\brho}{\boldsymbol{\rho}}
\newcommand{\bmu}{\boldsymbol{\mu}}

\begin{document}

\title{Some Results On Normal Homogeneous Ideals}
\author{Les Reid}
\address{Department of Mathematics\\Southwest Missouri State 
University\\ Springfield, MO
65804}
\email{les@math.smsu.edu}
\author{Leslie G. Roberts}
\address{Department of Mathematics and Statistics\\Queen's 
University\\ Kingston, Ontario,
CANADA K7L 3N6}
\email{robertsl@mast.queensu.ca}
\author{Marie A. Vitulli}
\address{Department of Mathematics\\University of Oregon, Eugene, OR 
97403}
\email{vitulli@math.uoregon.edu}
\subjclass{Primary 13C13; Secondary 13A20, 13F20}

\begin{abstract}
In this article we investigate when a homogeneous ideal in a
graded ring is normal, that is, when  all positive powers of the ideal
are integrally closed.  We are particularly interested in homogeneous 
ideals
in an $\N$-graded ring $A$ of the form $A_{ \ge m} := \bigoplus_{\ell 
\ge m}
A_{\ell}$ and monomial ideals in a polynomial ring over a field. For 
ideals of the 
form $A_{ \ge m}$ we generalize a recent result of Faridi.  We prove 
that a monomial 
ideal in a polynomial ring in $n$ indeterminates over a field is 
normal if and 
only if the first $n - 1$ positive powers of the ideal are integrally 
closed.  We
then specialize to the case of ideals of the form $I(\blambda) := 
\overline{J(\blambda)}$, where
$J(\blambda) = (x_1^{\lam_1}, \ldots, x_n^{\lam_n}) \subseteq 
K[\xvec{n}]$.  To state our main
result in this setting,  we let $\ell =  
\lcm(\lambda_1,\ldots,\widehat{\lambda_{i}},\ldots
\lambda_n)$, for $1 \le i \le n$, and set  
$\blambda^{\prime}=(\lambda_1,\ldots,\lambda_{i-1},
\lambda_i+\ell,\lambda_{i+1},\ldots,\lambda_n)$.  We prove that if 
$I(\blambda')$ is normal then
$I(\blambda)$ is normal and that the converse holds with a small 
additional
assumption. 
\end{abstract}
\maketitle

\section{Introduction}
In this paper we explore when a homogeneous ideal in a
graded ring is normal, that is, when all positive powers of the ideal 
are integrally closed.  In particular, we are interested in homogeneous 
ideals of an
$\N$-graded ring $A$ of the form $A_{ \ge m} := \bigoplus_{\ell \ge m}
A_{\ell}$ and monomial ideals in a polynomial ring over a field.  In 
the first
setting, we generalize a recent theorem of Faridi \cite{F}.  As for 
monomial
ideals, our first new result is that a monomial ideal $I$ in
a polynomial ring $K[\xvec{n}]$ over a field $K$ is normal if and 
only if the first
$n-1$ positive powers of $I$ are normal. We then specialize to the 
case of monomials ideals of the form
$\overline{J(\blam)}$, where
$J(\blam) := (x_1^{\lam_1}, \ldots, x_n^{\lam_n})$ is an ideal
in $R := K[\xvec{n}]$, $\blam := (\lam_1,
\ldots ,\lam_n)$ is a vector of positive integers, and  
$\overline{J(\blam)}$ is
the integral closure of ${J(\blam)}$ in $R$. 

In \cite{BG} Bruns and Gubeladze studied the normality of the
polytopal semigroup ring
$K[S(\blam)]$, where $K$ is a field and $S(\blam)$ is the submonoid of
$\N^{n+1}$ generated by 
$$\{ (a_1,\ldots,a_n,d) \in \mathbb{N}^{n+1} \mid a_1/
\lambda_1 + \cdots + a_n/ \lambda_n \le d \mbox{ for } d \le 1 \}.$$

Bruns and Gubeladze defined $\blam$ to be normal provided that
$K[S(\blam)]$ is normal.  One striking result in
\cite{BG} is the following theorem.

\begin{theorem}\label{main} {\rm \cite[Theorem 1.6]{BG}} 
Let $\blam = (\lam_1,\ldots ,\lam_n)$ be a vector of positive 
integers and set
$\ell = 
\lcm(\lambda_1,\ldots,\widehat{\lambda_{i}},\ldots
\lambda_n)$.  Then $\blam$ is normal if and only if   
$\blambda^{\prime}=(\lambda_1,\ldots,\lambda_{i-1},
\lambda_i+\ell,\lambda_{i+1},\ldots,\lambda_n)$ is normal; in other 
words the normality of $\blam$
depends only on the residue class of $\lam_i$ modulo the least common
multiple of the $\lam_j$ with $i \ne j$. 
\end{theorem}

Notice that Theorem \ref{main} says that the semigroup ring 
$K[S(\blam)]$ is
normal if and only if $K[S(\blam')]$ is normal.

The normality of the ideal $I(\blam) :=\overline{J(\blam)}$ is 
equivalent to
the normality of the semigroup ring $K[S'(\blam)]$, where $S'(\blam)$ 
is the
submonoid of
$\N^{n+1}$ generated by
$$ \{ (a_1,\ldots,a_n,d) \in \mathbb{N}^{n+1} \mid
a_1/
\lambda_1 + \cdots + a_n/ \lambda_n \ge d \mbox{ for } d \le 1 \}.$$
Due to the similarity between the semigroups
$S(\blam)$ and
$S^\prime(\blam)$ one might ask the following questions.

\begin{question}\label{conjecture0} Is $K[S(\blambda)]$ normal if and
only if $K[S^\prime(\blambda)]$ is normal (that is, if and only if
$I(\blam)$ is normal as an ideal)? 
\end{question}

\begin{question}\label{conjecture} Is $I(\blambda)$ normal if and
only if $I(\blambda^\prime)$ is normal? 
\end{question}

The answer to both of the above questions is no
by Example \ref{example4}.
Since we are interested in the normality of the ideal $I(\blam)$ 
rather than
the normality of the polytopal semigroup ring
$K[S(\blam)]$ and the normality of one does not imply the normality 
of the
other, we will no longer refer to the normality of the vector
$\blam$.  Later in the paper we identify the semigroup ring 
$K[S'(\blam)]$ with
the Rees algebra $R[I(\blam)t]$ and drop further references to 
$K[S'(\blam)]$. 
 The
normality of
$I(\blambda)$ for specific $\blambda$ can be  determined
readily using the {\bf normaliz} program \cite{normaliz} of Bruns and
Koch.

We now describe the organization of this paper.  In section 2 we 
review
some background material for our work, including
integral closure of monomial algebras and ideals, normality of ideals,
and polytopal semigroup rings. In section 3 we prove several results
on normal ideals in polynomial rings and $\N$-graded rings  of the
form $A_{ \ge m}$, generalizing
recent results of Faridi \cite{F}. In section 4 we develop for
$I(\blambda)$ an analogue of \cite[Proposition 1.3]{BG}. We introduce
the concept of quasinormality for an additive semigroup of the
nonnegative rational numbers. We  show that if
$I(\blam)$ is normal then the semigroup $\Lambda$ of $\Q_{\ge}$ 
generated
by $1/\lam_1, \ldots, 1/\lam_n$ is quasinormal (Lemma
\ref{ifnormalthenquasinormal}),
 and  if the $\lambda_i$ are pairwise relatively
prime then the converse holds
(Proposition \ref{normal iff qnormal}). In section 5 we show that the
two aforementioned questions have negative answers. Neither 
implication of
Question
\ref{conjecture0} holds.  This is shown in Example \ref{example4}.
However the implication  
$I(\blambda')$ normal implies  $I(\blambda)$ normal of Question  
\ref{conjecture} always  holds and the converse holds with an
additional hypothesis (Theorem \ref{congruence}). 
 
\smallskip
\textbf{Conventions.} All rings are assumed to be commutative with
identity.  
We let $\Z_+$ denote the set of positive integers, 
$\N$ the set of nonnegative integers,
$\Q_\ge$ the set of nonnegative rational numbers, $\Q_+$ the set of 
positive rational numbers, 
$\R_{\ge}$ the set of nonnegative real numbers, and $\mathbf{e}_1, 
\dots,
\mathbf{e}_n$ the standard basis vectors in $\R^n$.  We write 
$\balpha 
\le_{pr} \bbeta$ for vectors $\balpha = (a_{1},\ldots,a_{n}), \bbeta 
= 
(b_{1},\ldots,b_{n}) \in \R^{n}$ provided that $a_{i} \le b_{i}$ for 
 $1 \le i \le n$.  Thus $\balpha <_{pr} \bbeta$ means that $a_i \le 
b_i$ for all $1
\le i \le n$ and $a_j < b_j$ for some $1 \le j \le n$.  For a subset 
$X$ of $\R^n$ we
let $\conv(X)$ denote the convex hull of $X$.  Throughout this paper $R$ 
will 
denote the polynomial ring $K[x_1,\dots, x_n]$ 
over a field $K$ and $\blam = (\lam_1, \ldots, \lam_n)$ a vector of
positive integers.  In this context, for a vector 
$\balpha = (a_1, \ldots, a_n) \in \N^n$ we let $x^{\balpha}$ denote the monomial
$x_1^{a_1} \cdots x_n^{a_n}$.

\section{Background}

    In this section we recall some of the background to our 
investigation.
The integral closure of rings and ideals are as defined, for example,
in \cite[Chapter 4]{Eis}. We will be working primarily with monomial
ideals $I\subset R = K[x_1,\ldots,x_n]$ and subalgebras
$A\subset R$ that are generated by a finite number of 
monomials.  In these cases  we recall some definitions and notation that appeared in
\cite{RR} and \cite{RV}.
 
\begin{definition}\label{gamma} Let $X$ be any subset of 
$R=K[x_1,\ldots,x_n]$.
Then set
$$\Gamma(X) \; = \; \{\balpha \in \N{}^n \hspace{.1cm} \mid  
x^{\balpha} \in X\}.$$

\end{definition}

We refer to $\Gamma(X)$ as the {\em exponent set of $X$}. 
If $I$ is a monomial ideal then $\Gamma(I)$ is an ideal of the 
monoid $\mathbb{N}^n$ \cite[page 3]{Gil}. If $A$ is a subalgebra of
$R$ generated by monomials then $\Gamma(A)$ is a submonoid of
$\mathbb{N}^n$, and $A$ is isomorphic to the monoid ring 
$K[\Gamma(A)]$. 

\begin{definition} 
 For an arbitrary subset $\Lambda$ of  $\R{}^n$ and a positive 
 integer $m$ we let  
\begin{eqnarray*}
m \cdot \Lambda &=&  \{\lambda_1 + \cdots + \lambda_m \hspace{.1cm} 
\mid  
\lambda_i \in
\Lambda \hspace{.2cm} (i=1, \dots , m) \};  \hspace{.2 cm} 
\mathrm{and} \\  
m\Lambda &=&  \{ m\lambda \mid \lambda \in \Lambda \}.
\end{eqnarray*}

If $\Lambda = \Gamma(I)$ (respectively $\Gamma(A)$) then 
$\conv(\Lambda)$ will be denoted NP$(I)$ (respectively NP$(A)$),
and will be referred to as the {\em Newton polyhedron of $I$} 
(respectively, of $A$). 
\end{definition}

The integral closures of monomial ideals and subalgebras now have 
the following geometric descriptions.

\begin{theorem}\label{normalization}
{\rm (a) } Let $I$ be a monomial ideal in $R=K[x_1,\ldots,x_n]$. Then 
the
integral closure $\overline{I}$ of $I$ in  $R$ is 
the ideal defined by
$\Gamma(\overline{I})=\NP(I)\cap\mathbb{N}^n$ (so that 
$\NP(I)=\NP(\overline {I}))$.  Furthermore $$\Gamma(\overline{I})=
\{ \balpha \in \N{}^n |   m\balpha \in m \cdot \Gamma(I) \mbox{ for 
some }
 m \ge 1\}.$$

{\rm (b)}  Let $A$ be a subalgebra of $R$ generated by
a finite number of monomials. Then the
integral closure $\overline{A}$ of $A$ in  $R$ is 
the monoid ring defined by \newline
$\Gamma(\overline{A})=\NP(A)\cap\mathbb{N}^n$.  Furthermore 
$\NP(A)$ is the cone spanned
by $\Gamma(A)$ (or by the exponents of a (finite) set of algebra 
generators of $A$) and $$\Gamma(\overline{A})=
\{ \balpha \in \N{}^n |   m\balpha \in \Gamma(A) \mbox{ for some }
 m \ge 1\}.$$
\end{theorem} 
 
\begin{proof}
{\rm (a)} See Exercises 4.22, 4.23 in \cite{Eis}.

{\rm (b)} This is \cite[Proposition 6.1.2]{BH}. See also
\cite[3.1]{RR} for a form closer to what we want here.
\end{proof}

Polytopal semigroup rings, introduced in
\cite{BGT}, are examples of such monomial algebras. 
A polytopal semigroup ring is a monoid algebra
$K[S_P]$, where
$P$ is a polytope in $\mathbb{R}^n$ whose vertices have
integer coordinates and 
$S_P$ is the submonoid of $\mathbb{R}^{n+1}$
generated by the points $\{(\balpha,1) \mid \balpha \in P\cap\mathbb{Z}^n\}$.  The 
polytopal
semigroup ring of Bruns and Gubeladze that we referred to as 
$K[S(\blam)]$
in the Introduction is denoted in \cite{BG} by $K[S_{\Delta(\blam)}]$,
where the polytope
$\Delta(\blam) 
\subseteq
\R^n$ has vertices
$(0,\ldots,0), (\lam_1, 0, \ldots, 0), \ldots, (0, \ldots, 0, 
\lam_n)$.

An ideal $I$ in an integral domain $A$ is defined to be normal if
$I^m$ is integrally closed for all $m\in \mathbb{Z}_+$. The following
result is well known (for example, see \cite{Rib}).

\begin{theorem}\label{Rees} Let $I$ be an ideal in an integral domain 
$A$. 
Then the integral closure of the Rees algebra $A[It]$ in $A[t]$ is 
$\oplus_{i\geq 0}\overline{I^i}t^i.$ 
\end{theorem}

Thus in case the containing ring is a normal integral domain, the 
normality of $I$ is
equivalent to the normality of the Rees ring $A[It]$. Note that
``normal'' and ``integrally closed (in its quotient field)'' are 
synonyms for reduced
Noetherian rings but not for ideals. The following observation may be
helpful when contemplating normal ideals in $R=K[x_1,\ldots,x_n]$.

\begin{lemma} \label{lemles}Let $I = (x^{\bbeta_{1}},\ldots,
x^{\bbeta_{k}}) 
\subseteq R$ be a monomial ideal, let $m \ge 1$, and 
$J=(x^{m\bbeta_{1}},\ldots,x^{m\bbeta_{k}})$. Then
\begin{enumerate}
\item[{\rm (a)}]\hspace{1mm} $\NP(J)=\NP(I^{m}) =m\mathrm{NP}(I)= 
m\cdot\mathrm{NP}(I) = \newline \{ \balpha \in \R_{\ge}^{n} \mid 
\balpha  \ge_{pr} \sum_{i=1}^{k} c_{i}\bbeta_{i} {\mbox{ for some }}
 c_{i} \in \mathbb{R}_\ge, \sum c_{i} = m \}.$

\item[{\rm (b)}]\hspace{.1cm}  If  $\balpha \in \NP(I^{m})$ there 
exist $\, r$ affinely independent vectors $\bbeta_{i(1)},\ldots, 
\bbeta_{i(r)}$ in $\{ \bbeta_1,\ldots, \bbeta_k \}$ 
($r\leq n$) such that  $\balpha \in
m\mathrm{conv}(\bbeta_{i(1)},\ldots,\bbeta_{i(r)}) + \,  
\R_{\ge}^{n}$.

\end{enumerate}

\end{lemma}

\begin{proof}  (a) Obviously 
$\NP(J)\subseteq
\{\balpha\in\R_{\ge}^{n}\mid 
\balpha\ge_{pr}\sum_{i=1}^{k}c_{i}\bbeta_{i},
{\mbox{ for some }} c_{i}\ge 0,\sum c_{i} = m \}$ 
and the other sets mentioned in the lemma lie in between so it 
suffices to 
prove that  
$\{ \balpha \in \R_{\ge}^{n} \mid \balpha
\ge_{pr}\sum_{i=1}^{k}c_{i}\bbeta_{i}, c_{i}\ge 0,\sum c_{i} = m 
\}\subseteq
\NP(J)$. If 
$\bbeta=\sum_{i=1}^kc_i\bbeta_i + \bgamma$ with $c_i\geq 0$, $\sum 
c_i=m$ 
and $\bgamma \in \mathbb{R}_{\geq}^n$ then 
$\bbeta=\sum_{i=1}^k(c_i/m)(m\bbeta_i) + \bgamma$ with 
$\sum_{i=1}^kc_i/m=1$,
$m\bbeta_i \in \NP(J)$ and $\bgamma \in \mathbb{R}_{\geq}^n$ so 
$\bbeta \in \NP(J)$.

(b) If $\balpha\in\NP(I^m)$ then by (a) $\balpha=
\sum_{i=1}^kc_i\bbeta_i+\bgamma=m\sum_{i=1}^k(c_i/m)\bbeta_i+\bgamma$ 
where
$c_i\geq 0,\ \sum c_i=m$ (so that $\sum (c_i/m)=1$) and 
$\bgamma\in\mathbb{R}_\geq^n$. 
 By Carath\'{e}odory's Theorem $\bdelta:=\sum_{i=1}^k(c_i/m)\bbeta_i$ 
is in
the interior of \linebreak $\conv(\bbeta_{i(1)},\ldots,\bbeta_{i(r)})$
 where 
$\{\bbeta_{i(1)},\ldots,\bbeta_{i(r)}\}$
is an affinely independent subset of $\{\bbeta_{1},\ldots, \bbeta_{k} 
\}$
(so that $r\leq n+1$). If $r\leq n$ we are done. Otherwise
$r=n+1$  
 and there exists $t > 0$ such that 
$\bdelta - t\e_{1} 
\in  \, \conv(\bbeta_{i(1)},\ldots, \widehat{\bbeta_{i(j)}},\ldots, 
\bbeta_{i(n+1)})$ for some $1 \le j \le n+1$.  
Thus we may write $m\bdelta = 
\brho + \bnu$, where  $\brho \in m \, \conv(\bbeta_{i(1)},\ldots,
\widehat{\bbeta_{i(j)}},\ldots, 
\bbeta_{i(n+1)})$ and $\bnu \in \R_{\ge}^n$. Thus $\balpha=
\brho+(\bnu+\bgamma)
\in 
m\,\mathrm{conv}(\bbeta_{i(1)},\ldots,\widehat{\bbeta_{i(j)}},\ldots,
\bbeta_{i(n+1)}) + \,  \R_{\ge}^{n}.$
\end{proof}
\begin{remark} \label{normRees} It follows from the above discussion 
that 
the normalization  
$\overline{R[I(\blambda)t]}$ of $R[I(\blambda)t]$ is the subalgebra 
of 
$R[t]=K[x_1,\ldots,x_n,t]$ 
generated  by all $x^{\balpha} t^d$ where $\balpha=(a_1,\ldots,a_n)$, 
 $a_i, d\in\mathbb{N}$, and
\[ \frac{a_1}{\lambda_1} +\cdots+\frac{a_n}{\lambda_n} \geq  d.\]
On the other hand $\overline{K[S(\blam)]}$ is isomorphic to the 
subalgebra of 
$K[x_1,\ldots,x_n,t]$ 
generated  by all $x^{\balpha} t^d$ where $\balpha=(a_1,\ldots,a_n)$, 
 $a_i, d\in\mathbb{N}$, and
\[ \frac{a_1}{\lambda_1} +\cdots+\frac{a_n}{\lambda_n} \leq  d.\]
A crucial difference between the two cases is that 
$\overline{R[I(\blambda)t]}$ contains $x_1,\ldots,x_n$ but not $t$ 
whereas
$\overline{K[S(\blam)]}$ contains $t$ but not $x_1,\ldots,x_n$. 
\end{remark}
\medskip

\section{First Results on Normal Ideals}
In this section we obtain our first new result, namely that  a 
monomial ideal
$I$ in $R=K[x_1,\ldots,x_n]$ is normal if and only if the powers 
$I^m$ for 
$1\leq m\leq n-1$ are integrally closed. For the case $n = 2$ this 
follows from
the celebrated theorem of Zariski \cite[Appendix 5]{Zar} that asserts 
that the
product of integrally  closed ideals in a 2-dimensional regular ring 
is again
integrally closed.  Then we obtain a number of results
on the normality of the ideal
$A_{\ge d}$ of all elements of degree
$\geq d$ in the
$\mathbb{N}$-graded ring $A$. Note that the ideal 
$I(\blambda)$ is of this form for a suitable grading on 
$R$.

\begin{prop} \label{proples} Let $I \subseteq R=K[\xvec{n}]$ be a 
monomial 
ideal. If $I^{m}$ is integrally closed for $m = 1,\ldots,n-1$, then 
$I$ is 
normal.
\end{prop}

\begin{proof}   By Lemma \ref{normalization} part(a) it suffices to 
show
that if $\Gamma(I^{m}) = \mathrm{NP}(I^{m}) 
\cap \N^{n}$ for $m=1,\ldots,n-1$, then $\Gamma(I^{m}) = 
\mathrm{NP}(I^{m}) 
\cap \N^{n}$ for all $m \ge 1$. 

  Let $m \ge n$ and assume that $\Gamma(I^{i}) = \mathrm{NP}(I^{i}) 
\cap \N^{n}$ for all $i,\ 1\leq i\leq m-1$.  Clearly $\Gamma(I^{m}) 
\subseteq  \mathrm{NP}(I^{m}) 
\cap \N^{n}$.  Suppose that $\bgamma \in \mathrm{NP}(I^{m}) 
\cap \N^{n}$.  By Lemma  \ref{lemles} part (b),
$\bgamma \in m \, \mathrm{conv}(\bbeta_{1},\ldots,\bbeta_{r}) + 
\R^{n}_{\ge} \subseteq \mathrm{NP}(I)$ where 
$\bbeta_{1},\ldots,\bbeta_{r} \in \Gamma(I)$ are affinely 
independent  
and we may assume that $r \le n$.   Write $\bgamma = \sum_{i=1}^{r} 
c_{i}\bbeta_{i}  + \bnu$, where each $c_{i} \ge 0$, $\sum c_{i} = m$, 
and $\bnu \ge_{pr} 0$.  Since $m \ge n \ge r$, some $c_{i} \ge 1$.  
We may and shall assume that 
$c_{1} \ge 1$.  Then, $\bgamma - \bbeta_{1} = \sum_{i=1}^{r} 
(c_{i}-\delta_{i 1})\bbeta_{i} + \bnu $ and $\sum_{i=1}^{r} 
(c_{i}-\delta_{i 1}) = m -1$.  By the induction hypothesis, 
$\bgamma - \bbeta_{1} \in \Gamma(I^{m-1})$ and hence $\bgamma \in 
\Gamma(I^{m})$.  
\end{proof}

\begin{remark}The special case of Proposition 3.1 when $I$ is integral
over the  subideal generated by all monomials of the least total degree follows from
 ~\cite[Theorem 3.3]{BBV}.



\end{remark}
Another useful observation is the following.
\begin{lemma} \label{J=aI}  Let $A$ be a normal 
integral domain, $I \subseteq A$ be an ideal, and let $J = aI$, for 
$a \in A$.  Then, the following hold.
\begin{enumerate}
\item[(a)] $I$ is integrally closed if and only if $J$ is integrally 
closed.
\item[(b)] $I$ is normal if and only if $J$ is normal.
\end{enumerate}
\end{lemma}

\begin{proof}  Notice that an element $x \in A$ is integral over 
$J$ if and only if $x/a \in A$ and is integral over $I$. Part (a) 
and part (b) follow immediately.  \hspace{.1cm} 
\end{proof}

\begin{notn}  For an $\N$-graded ring $A$ and a positive integer $m$, 
we let $A_{\ge
m}$ denote  the homogeneous ideal defined by $A_{\ge m} = 
\bigoplus_{\ell \ge m} 
A_{\ell}$.
\end{notn}

Lemma \ref{faridilem} and Proposition \ref{faridithm} below 
generalize recent results of 
Faridi ~\cite[Lemma and Theorem 3]{F}. 

\begin{lemma} \label{faridilem} Let $A$ be an $\N$-graded ring 
generated 
over $A_{0}$ by homogeneous elements $x_{1}, \ldots,x_{n}$ of  
positive 
degrees $\omega_{1}, \ldots,\omega_{n}$ and $w = 
\lcm(\omega_{1},\ldots,\omega_{n})$.   Consider the ideal $I = A_{\ge 
kw}$ for a positive integer $k $.  
If $I^{p} = A_{\ge pkw}$ for $1 \le p \le \frac{n-2}{k} + 1$, then
$I^{p} =  A_{\ge pkw}$ for all $p \ge 1$.  In particular, if $k 
\ge n-1$, then $I^{p} =  A_{\ge pkw}$ for all $p \ge 1$.
\end{lemma}

\begin{proof}   We proceed by induction on $p \ge 1$, 
the case $p \le \frac{n-2}{k} + 1$ being a priori true.  

Suppose that $p > \frac{n-2}{k} + 1$ and that $I^{p-1} = 
A_{\ge (p-1)kw}$.  Let $\bmu = x_{1}^{c_{1}}\ldots x_{n}^{c_{n}}$ be 
a 
monomial of degree at least $pkw$. We must show that $\bmu \in I^{p}$.

Set $\lam_{i}=w/\omega_{i} \; (i=1,\ldots,n)$ and let
$q_i = \lfloor{c_{i} /\lam_{i}} \rfloor$  $(i=1,\ldots,n)$.    
Then
\begin{eqnarray*}
 pkw \le \deg(\bmu) &=& \sum_{i=1}^{n} c_{i}\omega_{i}  \\
 & < & \sum _{i=1}^{n} ( q_{i}+1)\lam_{i}\omega_{i}   \\
    & = & \left (\sum _{i=1}^{n} q_{i}\right )w + nw.   
\end{eqnarray*}
This implies that $ \sum q_{i} \ge pk -n +1$.   Since we assumed that 
$p > 
\frac{n-2}{k} + 1,$ we have  $(p-1)k \ge n-1$ and hence $pk 
-n +1 \ge k$.  Thus $\sum q_{i} \ge k$.

 Choose integers $0 
\le s_{i} \le q_{i} \; (i=1,\ldots,n)$ with $s_{1} + \cdots + s_{n} = 
k$.  Then, $\bmu = (x_{1}^{s_{1}\lam_{1}}\ldots
x_{n}^{s_{n}\lam_{n}})\bnu$  where $\deg\bnu =$  deg($\bmu) - (s_{1}+
\cdots + s_{n})w = \deg(\bmu) -kw \ge (p-1)kw$.  Thus, $\bnu \in 
I^{p-1}$
by the  induction hypothesis and 
$\bmu \in I^{p}$, as desired.  Since the inclusion $I^{p} \subseteq 
A_{\ge pkw}$ is immediate, the assertion is proven. 
\end{proof}

We suspect that the following result is well known but we do not know 
a 
reference so we provide a brief proof. We point out that this result 
holds if $\N$
is replaced by any totally ordered abelian group $G$ and we assume 
that $A$ is
positively graded.  The proof goes through without any changes.

\begin{lemma}\label{intclosed} Let $A$ be a reduced $\N$-graded ring 
and let $I = 
A_{\ge d}$ for some positive integer $d$.  Then, $I$ is an 
integrally closed ideal.
\end{lemma}

\begin{proof} Assume that $x \in A$ is integral over $I$ so that
$x^{n}+a_{1}x^{n-1}+ \cdots + a_{n} = 0$, for some $n \ge 1$ and 
$a_{k} \in I^{k} 
\; (k=1,\ldots,n)$.  Just suppose that the smallest component $x(i)$ 
of $x$ has degree $i < d$.  Since $a_{k} \in I^{k} \subseteq A_{\ge 
kd}$, the smallest component of $a_{k}x^{n-k}$ has degree strictly 
greater than $ni$ for $k = 1,\ldots,n$.  Hence we must have 
$x(i)^{n} = 0$, contradicting the assumption that $A$ is reduced. 
Thus $x \in I$. 
\end{proof}

\begin{prop} \label{faridithm} Let $A$ be a reduced $\N$-graded ring 
generated over $A_{0}$ by 
homogeneous elements $x_{1}, \ldots,x_{n}$ of  positive degrees 
$\omega_{1},
\ldots,\omega_{n}$ and $w = \linebreak
\lcm(\omega_{1},\ldots,\omega_{n})$.   Consider the ideal $I = A_{\ge 
kw}$ for a positive integer $k $.  
If $I^{p} = A_{\ge pkw}$ for $1 \le p \le \frac{n-2}{k} + 1$, then
 $I$ is a normal ideal.  In particular, if $k \ge n-1$, then 
$I$ is a normal ideal.  In this case, if $A$ is a normal 
domain, then the Rees ring $A[It]$ is again
a normal domain.

\end{prop}

\begin{proof}  This follows immediately from Lemmas 
\ref{faridilem} and \ref{intclosed}. 
\end{proof}

\begin{remark}   Notice that for $I = A_{\ge kw} \subseteq A$ as in 
Lemma \ref{faridilem}, we always have $\overline{I^{p}} = A_{\ge 
pkw}$.  
The containment $\overline{I^{p}} \subseteq A_{\ge pkw}$ follows 
from Lemma \ref{intclosed}.  To see the opposite containment, suppose 
that $x^{\bgamma} \in A_{\ge pkw}$, where 
$\bgamma=(c_1,\ldots,c_n)$.  
Then, $(x^{\bgamma})^{kw}= \prod x_{i}^{k\lam_{i}\omega_{i}c_{i}} = 
\prod ( x_{i}^{k\lam_{i}})^{\omega_{i}c_{i}} \in 
I^{\omega_{1}c_{1}}\cdots I^{\omega_{n}c_{n}} = I^{\bomega \cdot 
\bgamma}$.  However $x^{\bgamma} \in A_{\ge pkw}$ implies $\bomega \cdot 
\bgamma \ge pkw$ so that $(x^{\bgamma})^{kw} \in (I^{p})^{kw}$
and we are done.
\end{remark}

\section{$\mathfrak{m}$-Primary Monomial Ideals}

\textbf{Conventions.}  Let $\mathfrak{m}=(\xvec{n})$ denote the 
maximal homogeneous ideal of $R=K[x_1,\ldots,x_n]$. 
Furthermore let $\blam=(\lam_{1},\ldots,\lam_{n})$, 
$J(\blam) =(x_{1}^{\lam_{1}},\ldots,x_{n}^{\lam_{n}})$, and 
$I(\blam)=\overline{J(\blam)}$, as in the Introduction.
  We ask when 
the integrally closed $\mathfrak{m}$-primary monomial ideal 
$I(\blam)$ of $R$ is normal.

\begin{notn} \label{notation4}
 Let $L= \lcm(\lam_{1},\ldots,\lam_{n})$, 
$\omega_{i}=L/\lam_{i}$, 
$1/\blambda=(1/\lambda_1,\ldots,1/\lambda_n)$ 
and $\bomega = 
(\omega_{1},\ldots,\omega_{n})$, so that $L/\blambda=\bomega$. 
We will denote $\Gamma(I(\blam))$ 
(Definition \ref{gamma}) simply by $\Gamma$.
\end{notn}

Observe that 
$\NP(I(\blam))=\NP(J(\blam))$ has one bounded facet 
with vertices $\lam_{1}\e_{1},\ldots,\lam_{n}\e_{n}$. 
For $\balpha=(a_1,\ldots,a_n)$ the hyperplane 
$ (1/\blambda) \cdot \balpha =a_1/\lambda_1+\cdots+a_n/\lambda_n=1$
passes through these vertices, and upon multiplication by $L$, the
equation of  this hyperplane becomes $\bomega\cdot\balpha=L$. This
explains the lemma below.

\begin{lemma} \label{notation4a} 
$\Gamma = \{ \balpha\in \N^n \mid  (1/\blambda) \cdot \balpha \ge 1\}
= 
\{ \balpha \in \N^{n} \mid \bomega \cdot \balpha
\ge L\}$. 
\end{lemma}
Assigning deg($x_{i}) = \omega_{i} \; (i=1,\ldots, n)$ we now 
have $I(\blam) = R_{\ge \, L}$. Furthermore, the   
following gives  necessary and sufficient conditions for 
$I(\blam)$ to be normal.

\begin{lemma} \label{lem reductions} For the ideal $I(\blam)$ defined 
above the  following are equivalent.  
\begin{enumerate}
\item[{\rm (a)}] $I(\blam)$ is normal. 
\item[{\rm (b)}]   Whenever $\bomega \cdot \balpha \ge pL$ for 
$\balpha 
\in \N^{n}$ and $p \in \N$, there exist vectors \\ $\bbeta_{j} \in 
\Gamma \; 
(j=1,\ldots,p)$ 
 such that   $\balpha = \sum \bbeta_{j}.$
 \item[(c)]  Whenever $\bomega \cdot \balpha \ge pL$ for 
$\balpha=(a_1,\ldots,a_n)\in\mathbb{N}^n$ with $\lambda_i>a_i$ and 
\linebreak
$1\leq p<n$, there exist vectors $\bbeta_{j}
\in
\Gamma
\; (j=1,\ldots,p)$ 
 such that   $\balpha = \sum \bbeta_{j}.$
  
\end{enumerate}
\end{lemma}

\begin{proof} The equivalence of (a) and (b) is an immediate 
consequence 
of Lemma \ref{lemles} and Theorem \ref{normalization}. Clearly (b) 
implies (c) so it remains only to show that (c) implies (b).

Suppose (c) holds. 
We need only verify  condition (b) for $2\leq p < n$ 
by Lemma \ref{proples} and the observation that (b) automatically 
holds for 
$p=1$.  We argue by decreasing induction on $p$.
Assume that $\bomega\cdot\balpha\geq pL$ for $\balpha \in \N^{n}$.
If $a_i < \lambda_i$ for all $i$ we can apply (c) directly. If
 $a_{i} \ge \lam_{i}$ for some $i \in \{1,\ldots,n \}$ (we may assume 
that $a_{1} \ge \lam_{1}$) then $\balpha = (\balpha - \lam_{1}\e_{1}) 
+ 
\lam_{1}\e_{1}$. Dotting with $\bomega$ we obtain 
$\bomega \cdot \left( (\balpha - \lam_{1}\e_{1}) + 
\lam_{1}\e_{1} \right) \ge pL$, which implies 
 $\bomega \cdot  (\balpha - 
\lam_{1}\e_{1})  \ge (p-1)L$ (since $\bomega\cdot \lambda_1\e_1=L$).  
By induction, there exist vectors 
$ \bbeta_{j} \in 
\Gamma\ (j=1,\ldots, p-1)$ with $\balpha - \lam_{1}\e_{1} = \sum 
\bbeta_{j}$.  Thus $\balpha = \sum \bbeta_{j} + \lam_{1}\e_{1}$ and 
condition (b) is satisfied. 
\end{proof}

\medskip

  Due to this characterization, we will say that $\Gamma$ is 
  \emph{normal} if either condition (b) or (c) above holds (so that
$\Gamma$ is normal if and only if $I(\blambda)$ is normal).

To put this section into context with 
the preceding section notice that if $w = 
\lcm(\omega_{1},\ldots,\omega_{n})$ and 
$d=\gcd(\lam_{1},\ldots,\lam_{n})$ then $L = dw$ 
(this equality is easily checked by showing that any prime number $p$
has the same exponent in $L$ and in $dw$)
so that $I(\blam)= 
R_{\ge dw}$.  From this point of view we obtain the following 
corollary.

\begin{cor}  Let $\blam=(\lam_{1},\ldots,\lam_{n}) \in \Z_+^{n},\ 
n\geq 3,$ 
and suppose that \newline $\gcd(\lam_{1},\ldots,\lam_{n}) >n-2$.  
Then the 
monomial ideal $I(\blam) \subseteq K[x_{1},\ldots,x_{n}]$ is normal. 
\end{cor}
\begin{proof}  This is an immediate consequence 
of Proposition \ref{faridithm} since, in the notation of that result, 
$k \ge   n-1$. 
\end{proof}

In \cite{BG}  Bruns-Gubeladze define a  submonoid 
$S$ 
 of $\Q_{\ge}$ to be  \emph{1-normal} if  
whenever $x \in S$ and $x \leq p$ for some $p \in \N$, there exist 
rational numbers $y_{1},\ldots, y_{p}$ in $S$ with $y_i\leq 1$ for 
all $i$ 
such  that $x = y_{1} + \cdots + y_{p}$. Then they relate the 
normality of
$K[S(\blam)]$ to the 1-normality of the submonoid  
$\Lambda$ (defined below) of $\mathbb{Q}_{\geq}$.  
We modify this program as follows.

\begin{definition}
A submonoid  $S$  of $\Q_{\ge}$ is \emph{quasinormal} provided that 
whenever $x \in S$ and $x \ge p$ for some $p \in \N$, there exist 
rational numbers $y_{1},\ldots, y_{p}$ in $S$ with $y_i\geq 1$ for 
all $i$ 
such that $x = y_{1} + \cdots + y_{p}$.
\end{definition}
 
We now have the following.

\begin{lemma}\label{ifnormalthenquasinormal}

   Let 
$\Lambda = \langle 1/\lam_{1},\ldots,1/\lam_{n}\rangle$, the additive 
submonoid of $\Q_{\ge}$   generated by $1/\lam_{1},\ldots,1/\lam_{n} 
\,$,
and 
$\Lambda_{\ge 1} = \{ x\in \Lambda \mid x \ge 1 \}$.
If $I(\blambda)$ is normal then $\Lambda$ is quasinormal.

\end{lemma}
\begin{proof} 
Suppose $I(\blambda)$ is normal
and $x\in\Lambda,\ x\geq p$.
Then $x= (1/\blambda) \cdot  \balpha$ for $\balpha\in\mathbb{N}^n$. 
As noted in  Lemma \ref{notation4a},
$\bomega \cdot \balpha \geq pL$. Therefore by Lemma \ref{lem 
reductions}
there exist vectors $\bbeta_i\in\Gamma,\ 1\leq i\leq p$, so that 
$\balpha = \bbeta_{1} + \cdots + \bbeta_{p}$. Thus 
$x= (1/\blambda) \cdot \balpha = (1/\blambda) \cdot \bbeta_{1}  + 
\cdots + 
(1/\blambda) \cdot \bbeta_{p} $. Again by the description of $\Gamma$ 
in 
Lemma \ref{notation4a}, $(1/\blambda)\cdot  \bbeta_i \in
\Lambda_{\geq 1}$. Hence $\Lambda$ is quasinormal.
\end{proof}

   When we assume that the 
integers $\lam_{1},\ldots,\lam_{n}$ are pairwise relatively prime the 
converse is true.  So in this special case, the normality condition 
on the 
$n$-dimensional monoid $\Gamma$ is reduced to the quasinormality 
condition on the 1-dimensional monoid $\Lambda$.

\begin{prop} \label{normal iff qnormal} Suppose that $\lam_{1},\ldots,
\lam_{n}$ are
pairwise  relatively prime positive integers and let $\Lambda$ be as 
in
Lemma 
\ref{ifnormalthenquasinormal}.  
With notation and assumptions  as above, 
$I(\blam)$ is normal if and only if $\Lambda$ is quasinormal.  
\end{prop}

\begin{proof}  
By Lemma \ref{ifnormalthenquasinormal} it suffices to show that if 
$\Lambda$
is quasinormal then $I(\blambda)$ is normal. We assume that $\Lambda$
is quasinormal and establish the criterion 
for normality of $I(\blambda)$
in Lemma \ref{lem reductions}(c).

First, as in ~\cite[Proposition 1.3]{BG}, we consider the natural 
surjection 
$$\pi: \Z^{n} 
\rightarrow \mathrm{grp}(\Lambda),$$ defined by  $\e_{i} \mapsto 
1/\lam_{i} 
\; (i=1,\ldots,n)$, where $\mathrm{grp}(\Lambda)$ is the subgroup of 
$\mathbb{Q}$ generated by $\Lambda$. 
Suppose $(a_1,\ldots,a_n)\in \ker(\pi)$. Clearing denominators in the
equation $a_1/\lambda_1+\cdots+ a_n/\lambda_n=0$ and using that the
$\lambda_i$ are pairwise relatively prime we observe that $\lambda_i$
divides
$a_i$.   In particular any nonzero element of $\ker(\pi)$ has $i^{{\rm
th}}$ coordinate greater than or equal to $\lambda_i$ for some $i$. 
  
Suppose that $\balpha = (a_{1},\ldots,a_{n})  \in \N^{n}$ 
satisfies  $\bomega \cdot \balpha \ge pL$, with $\balpha <_{pr} 
\blam$, 
as in the hypotheses of Lemma \ref{lem  reductions}(c).   
We have 
  $$ x := (1/\blambda)  \cdot \balpha = a_{1}/\lam_{1} + \cdots +
 a_{n}/\lam_{n} \ge p.$$ Since $\Lambda$ is quasinormal there exist 
numbers 
$y_{i} \in \Lambda, y_i\ge 1 
 \, (i =1,\ldots,p)$ such that $x = y_{1} + \cdots + y_{p}$.  
 Write $y_{i} =  (1/\blam) \cdot \bbeta_{i}$ with $\bbeta_{i} \in
\mathbb{N}^n 
 \ (i=1,\ldots,p)$. By Lemma 4.2, $\bbeta_i\in \Gamma$ for all $i$. 
  Then $\bbeta := \bbeta_{1} + \cdots + \bbeta_{p} \in  \balpha + 
 \mathrm{ker}(\pi)$. By the above 
discussion of $\ker(\pi)$, $\balpha$ is 
the only element in its coset  with nonnegative coordinates. Hence 
$\balpha = \bbeta$, which completes the proof.  
  \end{proof}

 \begin{prop} \label{1+1/l in Lambda} Let $\blam \in \Z_+^{n}$.
 If the monoid $\Lambda = \langle 1/\lam_{1},\ldots, 
 1/\lam_{n}\rangle$ is quasinormal, then $1 + 1/L \in \Lambda$.  
 
 \end{prop}
 
 \begin{proof} Assume that $\Lambda$ is quasinormal. Notice that 
 $\omega_{1},\ldots,\omega_{n}$ are 
 relatively prime and 1 is the smallest positive integer in 
 grp$(\omega_{1},\ldots,\omega_{n})$.  Hence $1/L$ is the smallest 
 positive number in grp($\Lambda$).  Choose an integer $ N \gg 0$ in 
 $\Lambda$ such 
 that $N + 1/L\in \Lambda$.  Write $N + 1/L= 
 (1/\blam) \cdot \bbeta $, where $\bbeta \in \N^{n}$.  
 Since $\Lambda$ is quasinormal, $N + 1/L= y_{1} + \cdots + y_{N}$, 
 where $y_{i} \in \Lambda, \, y_i\geq 1, \, (i = 1,\ldots,N)$.
 
 We claim that if $x \in \Lambda$ and 
 $x > 1$, then $x \ge 1 + 1/L$.  
  Now suppose that $1 < x \in \Lambda$.  Since $1 = 
 \lam_{1}(1/\lam_{1}) \in \Lambda$, $x - 1 \in \grp(\Lambda)$. 
Furthermore 
 $x - 1 > 0$ implies $x-1 \ge 1/L$.  Thus $ y_{1} + \cdots + y_{N} = 
 N + 1/L$  forces $N-1$ of the  $y_{i}$ to be one 
 and the remaining to be $1+1/L$.  In particular, $1 + 1/L
 \in \Lambda$.
\end{proof}

 Multiplying by $L$ we obtain the following version of the corollary.

 \begin{prop}   Let $\blam \in \Z_+^{n}$.
If $\Lambda = \langle 1/\lam_{1},\ldots,
 1/\lam_{n} \rangle$ is quasinormal, then $L + 1 \in \langle 
 \omega_{1},\ldots,\omega_{n}\rangle$.
 \end{prop}
 
\section{Normality of $I(\blambda)$ and $I(\blambda^\prime)$}

In this section we discussion the relationship between the normality 
of
$I(\blambda)$ and that of $I(\blambda^\prime)$. 
Our notation continues as usual: $R=K[x_1,\ldots,x_n]$ for a field 
$K$,
 $\blambda=(\lambda_1,\ldots,\lambda_n)$ for arbitrary positive 
integers 
$\lambda_j$, $L=\lcm(\lam_1, \ldots, \lam_n)$,  and $\blam' = 
(\lam_1, \ldots, \lam_{i-1},
\lam_i + \ell, \lam_{i+1} \ldots,
\lam_n)$, where $ \ell = \lcm(\lam_1, \ldots, \widehat{\lam_i},
\ldots, \lam_n)$.  There is no loss of generality in taking
$i=n$, so that
$\ell=\lcm(\lambda_1,\ldots,\lambda_{n-1})$ and $\blambda^\prime= 
(\lambda_1,\ldots,\lambda_{n-1},\lambda_n+\ell)$.
We now state our main result.

\begin{theorem} \label{congruence} If $I(\blambda^{\prime})$ is 
normal then 
$I(\blambda)$ is normal. If $\lambda_n \geq \ell$ and  
$I(\blambda)$ is normal so is  $I(\blambda^{\prime})$.
\end{theorem}

Before beginning the proof of Theorem \ref{congruence} 
we will give an 
example to show that Question \ref{conjecture} has a 
negative answer.

 \begin{example}\label{example4}  The ideal $I = 
I(2,3,7) =
\overline{(x^{2},y^{3},z^{7})} \subseteq  K[x,y,z]$ is not 
normal.
 In this case, $L= 42$ and $L + 1 = 43$ is not in the 
 monoid generated by $\omega_{1}=21, \omega_{2}=14$, and 
 $\omega_{3}=6$.  Hence the monoid $\Lambda$ is not quasinormal, 
 which implies, by Lemma \ref{ifnormalthenquasinormal} that the ideal 
$I$ is not normal. Alternatively, 
$\balpha=(1,2,6)$ satisfies $(1/\blambda) \cdot \balpha \ge 2$ but
$\balpha$ is  not the sum of two elements of $\Gamma(I)$. Thus, by the
discussion of section 2, $xy^2z^6$ is integral over $I^2$ but not in
$I^2$. Hence $I^2$ is not integrally closed, so $I$ is not normal.
However by \cite[Theorem 1.6]{BG}, the ring $K[S(2,3,7)]$ is normal 
because  the ring $K[S(2,3,1)]$ is normal. 
Thus $\blambda=(2,3,7)$ is an example where $K[S(\blambda)]$ is normal
but
$I(\blambda)$ is not. Also $I(2,3,1)$ is normal, so 
we also have a  counterexample to Question \ref{conjecture0}.
In the other direction $K[S(2,3,5)]$ is not normal but $I(2,3,5)$ is
normal.
\end{example} 

\medskip
Our method of proof of Theorem \ref{congruence} is to compare the 
minimal
generators of the integral closures of the Rees algebras 
$R[I(\blambda)t]$
and  $R[I(\blambda^\prime)t]$. The integral closure  
$\overline{R[I(\blambda)t]}$ of $R[I(\blambda)t]$, by 
Remark \ref{normRees}, is the subalgebra of 
$R[t]=K[x_1,\ldots,x_n,t]$ 
generated  by all $x^{\balpha} t^d$ where $\balpha=(a_1,\ldots,a_n),$ 
such that $a_i, d\in\mathbb{N}$ and 

\medskip
\begin{enumerate}
\item [(A)] \hskip 1.5in 
${\displaystyle \frac{a_1}{\lambda_1} +\cdots+\frac{a_n}{\lambda_n} 
\geq  d.
} $
\end{enumerate}

\medskip
The algebra  $\overline{R[I(\blambda)t]}$ has a unique (finite) 
minimal  
set of monomial generators, 
corresponding to exponent vectors $(a_1,\cdots,a_n,d) \in \N^{n+1}$
of the following types:   

\begin{itemize}
\item []
\begin{enumerate}
\item The ``trivial'' exponent vectors  
   \[\begin{array}{llllll}1& 0& 0& \cdots& 0&0 \\
                       0& 1& 0& \cdots& 0&0 \\
                       &&&\cdots&& \\
                       0&0&0& \cdots & 1 &0 \end{array}\]
corresponding to $x_i\in R\subset \overline{R[I(\blambda)t]}$.
\item  
Exponent vectors of the form
 \[\begin{array}{llllll}\lambda_1& 0& 0& \cdots& 0&1 \\
                       0&\lambda_2& 0& \cdots& 0&1 \\
                       &&&\cdots&& \\
                       0&0&0& \cdots & \lambda_n &1 \end{array}\]
corresponding to $x_i^{\lambda_i}t\in I(\blambda)t\subset 
\overline{R[I(\blambda)t]}$.
\item Exponent vectors 
\[ (a_1,\ldots,a_{n-1},0,d) \]
  with 
 $d>0$ and $a_ia_j >  0$ for some $0 < i < j < n$.

\item Exponent vectors
\[ (a_1,\ldots,a_n,d) \] 
  with 
$d>0$  and $a_ia_n >0$ for some $0 < i < n$.

\end{enumerate}
\end{itemize}
The exponent vectors of type (2) have been written down separately 
because they are the initial data of the problem. 
The condition
that an exponent vector correspond to a minimal generator is that it 
cannot be written as the sum of two nonzero vectors satisfying 
condition (A).

     In the sequel we will
informally refer to the exponent vectors themselves as being 
generators of
$\overline{R[I(\blambda)t]}$. In this language 
generators of types (1) and (2) are obviously minimal.
The condition that $(a_1,\ldots,a_n,d)$ of type (3) or (4) be minimal 
is
that in addition to satisfying condition (A) it also satisfy:

\begin{enumerate}
\item[(B)] If any one of $a_1,...,a_n$ which is greater than 0 is
decreased, then inequality (A) fails to hold, and

\item[(C)] If $d>1$ then $(a_1,\ldots,a_n,d)$ cannot be written
in the form $(a_1,\ldots,a_n,d) = 
(b_1,\ldots,b_n,d_1)+(c_1,\ldots,c_n,d_2)$
with $0<d_1, d_2 <d$, where $(b_1,\ldots,b_n,d_1)$ and 
$(c_1,\ldots,c_n,d_2)$ both satisfy (A).
\end{enumerate}

Condition (B) says that $(a_1,\ldots,a_n,d)$ cannot be written as
the sum of a vector of type (1) and another vector satisfying (A).
Condition (C) says that $(a_1,\ldots,a_n,d)$ cannot be written as
the sum of two vectors of types (2), (3), or (4).  

In this context 
Theorem \ref{Rees} can be restated as follows.

\begin{lemma}\label{normalitycriterion}
$I(\blambda)$ is normal if 
and only if the minimal generators of types
{\rm (3)} and {\rm (4)} all have 
$d=1$.
\end{lemma}

The lemma below will be useful.

\begin{lemma}\label{finiterange} In any minimal generator of type
{\rm (3)}
or {\rm (4)} we have $0\leq a_i<\lambda_i$ for all $i$.
\end{lemma}

\begin{proof} We argue by contradiction. 
Let $(a_1,\ldots,a_n,d)$ be 
a minimal generator of type
(3) or (4) with $a_i\geq \lambda_i$ for some $i$. 
We may suppose without loss of generality that $a_1\geq\lambda_1$. 
From 
condition (A) we must have $d\geq 1$. If we subtract the equality
\[ \frac{\lambda_1}{\lambda_1} = 1\]
from the inequality (A) we obtain
\[ \frac{a_1-\lambda_1}{\lambda_1}+\frac{a_2}{\lambda_2}+\cdots + 
\frac
{a_n}{\lambda_n} \geq d-1,\]
that is,  $(a_1-\lambda_1,a_2, \ldots,a_n,d-1)$ satisfies condition 
(A) and
\newline
$(a_1-\lambda_1,a_2,\ldots,a_n,d-1)+(\lambda_1,0,\ldots,0,1)=(a_1,\ldots,a_n,d)$
so $(a_1,\ldots,a_n,d)$ is not minimal,  a contradiction.
\end{proof}

\begin{cor}\label{corrange} For any minimal generator 
$(a_1,\ldots,a_n,d)$ 
we have $d<n$.
\end{cor} 
Corollary \ref{corrange}  reproves Proposition \ref{proples}, but 
only in 
the special case of $I(\blambda)$.

We now compare the minimal generators of $\overline{R[I(\blambda)t]}$
and $\overline{R[I(\blambda^{\prime})t]}$. There is obviously a 
bijection
between minimal generators of types (1), (2), and (3) for these two 
algebras. If $a_n=0$ then $(a_1,\ldots,a_n,d)$ corresponds to itself,
as does $(0,0,\ldots,1,0)$, and the generator 
$(0,\ldots,0,\lambda_n,1)$
for $\overline{R[I(\blambda)t]}$ corresponds to the generator
$(0,\ldots,0,\lambda_n+\ell,1)$ for 
$\overline{R[I(\blambda^{\prime})t]}$.
Now let $(a_1,\ldots,a_n,d)$ be a minimal generator of 
$\overline{R[I(\blambda)t]}$ of type (4), i.e., with $a_n>0,d>0$. 
Then,
by condition (B), $a_n$ is the smallest integer such that (for
fixed $a_1,\ldots,a_{n-1},d$)
\[\frac{a_1}{\lambda_1}+\cdots+\frac{a_n}{\lambda_n}\geq d.\]   
Now define $a_n^{\prime}$ to be the smallest integer such that
\[\frac{a_1}{\lambda_1}+\cdots+
\frac{a_{n-1}}{\lambda_{n-1}}+\frac{a_n^{\prime}}{\lambda_n+\ell}\geq 
d,\]
so that $(a_1,\ldots,a_{n-1},a_n^{\prime},d)$ is the exponent vector 
of
an element of $\overline{R[I(\blambda^{\prime})t]}$. 

\begin{prop} \label{injection} Let $(a_1,\ldots,a_n,d)$ be a minimal 
generator of
$\overline{R[I(\blambda)t]}$ of type {\rm (4)}, and let $a_n^{\prime}$ be 
as 
defined above. Then $a_n^{\prime} =a_n+ 
d\ell-\frac{\ell}{\lambda_1}a_1-\cdots-\frac{\ell}{\lambda_{n-1}}a_{n-1}$ 
and
$(a_1,\ldots,a_{n-1},a_n^{\prime},d)$ is a minimal generator of
$\overline{R[I(\blambda^{\prime})t]}$ of type {\rm (4)}.
\end{prop}
\begin{proof}
Let $\delta$ be any integer, and consider the following chain of 
equivalent
inequalities:
\[\frac{a_1}{\lambda_1}+\cdots+\frac{a_{n-1}}{\lambda_{n-1}}+\frac{a_n+\delta}{\lambda_n+\ell}\geq 
d\Leftrightarrow\]
\[\ell\frac{a_1}{\lambda_1}+\cdots+\ell\frac{a_{n-1}}{\lambda_{n-1}}+\frac{\ell(a_n+\delta)}{\lambda_n+\ell}\geq 
d\ell\Leftrightarrow\]
\[\frac{\ell(a_n+\delta)}{\lambda_n+\ell}\geq 
d\ell-\frac{\ell}{\lambda_1}a_1-\cdots-\frac{\ell}{\lambda_{n-1}}a_{n-1}\Leftrightarrow\]
\[\ell a_n+\ell \delta\geq
(\lambda_n+\ell)(d\ell-\frac{\ell}{\lambda_1}a_1-\cdots-\frac{\ell}{\lambda_{n-1}}a_{n-1}).\]

In the rest of this proof we will set 
\[\delta = d\ell-\frac{\ell}{\lambda_1}a_1-
\cdots-\frac{\ell}{\lambda_{n-1}}a_{n-1}\] 
(which is an integer by the definition of $\ell$). Then the last 
inequality becomes
\[\ell a_n\geq \lambda_n(d\ell-\frac{\ell}{\lambda_1}a_1-\cdots-
\frac{\ell}{\lambda_{n-1}}a_{n-1}),\]
which is equivalent to 
\[\frac{a_1}{\lambda_1}+\cdots+\frac{a_{n-1}}{\lambda_{n-1}}+
\frac{a_n}{\lambda_n}\geq d.\]
Putting these equivalences together we conclude that for {\em any} 
integers
$a_1,\ldots,a_n$ and 
\[\delta = d\ell-\frac{\ell}{\lambda_1}a_1-
\cdots-\frac{\ell}{\lambda_{n-1}}a_{n-1}\]  
 we have
\flushleft $(*)$ \hskip .5in 
${\displaystyle 
\frac{a_1}{\lambda_1}+\cdots+\frac{a_{n-1}}{\lambda_{n-1}}+\frac{a_n+\delta}{\lambda_n+\ell}\geq 
d\Leftrightarrow 
\frac{a_1}{\lambda_1}+\cdots+\frac{a_{n-1}}{\lambda_{n-1}}+
\frac{a_n}{\lambda_n}\geq d.}$

\smallskip
If 
\[\frac{a_1}{\lambda_1}+\cdots+\frac{a_{n-1}}{\lambda_{n-1}}+
\frac{a_n+\delta-1}{\lambda_n+\ell}\geq d\]
then we have, by $(*)$ applied to
$a_1,\ldots,a_{n-1},a_n-1$, that
\[\frac{a_1}{\lambda_1}+\cdots+\frac{a_{n-1}}{\lambda_{n-1}}+
\frac{a_n-1}{\lambda_n}\geq d\]
which contradicts condition (B) for the minimality of  
$(a_1,\ldots,a_n,d)$ as a generator of
$\overline{R[I(\blambda)t]}$. Thus we have 
$a_n^{\prime} =a_n+ d\ell-\frac{\ell}{\lambda_1}a_1-\cdots
-\frac{\ell}{\lambda_{n-1}}a_{n-1}$ as claimed.
By construction $(a_1,\ldots,a_{n-1},a_n^{\prime},d)$ satisfies 
condition
(A) for minimality as a generator of 
$\overline{R[I(\blambda^{\prime})t]}$. If in 
$(a_1,\ldots,a_{n-1},a_n^{\prime},d)$ we replace $a_n^{\prime}$ by 
$a_n^{\prime}-1$ then condition (A) fails to hold (by the definition
of $a_n^{\prime}$). Because $a_n^{\prime}$ has the largest denominator
in the inequality (A), decreasing any of $a_1,\ldots,a_{n-1}$
will also violate condition (A). Therefore
$(a_1,\ldots,a_{n-1},a_n^{\prime},d)$ also satisfies condition (B)
for minimality as a generator of 
$\overline{R[I(\blambda^{\prime})t]}$.

Now we consider condition (C) for minimality of 
$(a_1,\ldots,a_{n-1},a_n^{\prime},d)$ as a generator of 
$\overline{R[I(\blambda^{\prime})t]}$. Note that 
since $(a_1,\ldots,a_n,d)$ is a minimal generator of 
$\overline{R[I(\blambda)t]}$ with $a_n>0$ we must have
\[\frac{a_1}{\lambda_1}+\cdots+\frac{a_{n-1}}{\lambda_{n-1}}< d\]
and hence that $\delta>0$. Furthermore $\delta$ is linear in
$a_1,\ldots,a_{n-1}$, and $d$. Hence we have an isomorphism of abelian
groups $f:\mathbb{Z}^{n+1}\to\mathbb{Z}^{n+1}$ defined by
\[f(u_1,\ldots,u_n,u_{n+1})=(u_1,\ldots,u_{n-1},
u_n+d\ell-\frac{\ell}{\lambda_1}u_1-\cdots-
\frac{\ell}{\lambda_{n-1}}u_{n-1},u_{n+1})\] which satisfies 
$f(a_1,\ldots,a_n,d)=(a_1,\ldots,a_{n-1},a_n^{\prime},d)$. 
Suppose that
 $(a_1,\ldots,a_{n-1},a_n^{\prime},d)$ fails to satisfy (C). Then we 
can write
\[(a_1,\ldots,a_{n-1},a_n^{\prime},d)=(b_1,\ldots,b_{n-1},b_n,d_1)+(c_1,\ldots,c_{n-1},c_n,d_2)\]

with $b_i,c_i\geq 0$, $0<d_1, d_2<d$ and 
\[\frac{b_1}{\lambda_1}+\cdots+\frac{b_{n-1}}{\lambda_{n-1}}+
\frac{b_n}{\lambda_n+\ell}\geq d_1,\]
\[\frac{c_1}{\lambda_1}+\cdots+\frac{c_{n-1}}{\lambda_{n-1}}+
\frac{c_n}{\lambda_n+\ell}\geq d_2.\]
Applying $f^{-1}$ we get
\[(a_1,\ldots,a_{n-1},a_n,d)=(b_1,\ldots,b_{n-1},b_n-\delta_1,d_1)+(c_1,\ldots,c_{n-1},c_n-\delta_2,d_2)\]

where 
\[\delta_1=d_1\ell-\frac{\ell}{\lambda_1}b_1-\cdots-\frac{\ell}{\lambda_{n-1}}b_{n-1}\]

and
\[\delta_2=d_2\ell-\frac{\ell}{\lambda_1}c_1-\cdots-\frac{\ell}{\lambda_{n-1}}c_{n-1}.\]

This will contradict the minimality of $(a_1,\ldots,a_n,d)$ as a 
generator of
$\overline{R[I(\blambda)t]}$ if we can show that $b_n-\delta_1\geq 0$ 
and
$c_n-\delta_2\geq 0$, and that
$(b_1,\ldots,b_{n-1},b_n-\delta_1,d_1)$, 
$(c_1,\ldots,c_{n-1},c_n-\delta_2,d_2)$
satisfy condition (A) for $\overline{R[I(\blambda)t]}$.
By $(*)$ we have
 \[\frac{b_1}{\lambda_1}+\cdots+\frac{b_{n-1}}{\lambda_{n-1}}+
\frac{b_n-\delta_1}{\lambda_n}\geq d_1\] and
 \[\frac{c_1}{\lambda_1}+\cdots+\frac{c_{n-1}}{\lambda_{n-1}}+
\frac{c_n-\delta_2}{\lambda_n}\geq d_2,\]
which is condition (A). 
 If $\delta_1\leq 0$ and $\delta_2\leq 0$ 
(which can happen) we will certainly have $b_n-\delta_1\geq 0$ and
$c_n-\delta_2\geq 0$.  Hence suppose that
$\delta_1>0$. Then 
 \[\frac{b_n}{\lambda_n+\ell}\geq 
d_1-\frac{b_1}{\lambda_1}-\cdots-\frac{b_{n-1}}{\lambda_{n-1}}=
\frac{\delta_1}{\ell},\]
so $b_n\geq (\lambda_n+\ell)\delta_1/\ell>\delta_1$ and 
$b_n-\delta_1> 0$.
Similarly if $\delta_2>0$ then $c_n-\delta_2> 0$.
This shows that $(a_1,\ldots,a_{n-1},a_n^{\prime},d)$ is a minimal
generator of $\overline{R[I(\blambda^{\prime})t]}$. Finally
$(a_1,\ldots,a_{n-1},a_n^{\prime},d)$ is of type (4) because 
$\delta>0$,
and hence a fortiori $a_n^{\prime}>0$.
\end{proof}
Now we show that if $\lambda_n\geq \ell$ then $f$ gives a bijection 
on 
minimal generators of type (4).
\begin{prop} \label{surjection} Let 
$(a_1,\ldots,a_{n-1},a_n^{\prime},d)$
be a minimal generator of
$\overline{R[I(\blambda^{\prime})t]}$ of type {\rm (4)}. Suppose that 
$\lambda_n\geq\ell$ and that $f$ is as defined in the proof of 
Proposition  \ref{injection}. Then
$f^{-1}(a_1,\ldots,a_{n-1},a_n^{\prime},d)$ is a minimal generator of
$\overline{R[I(\blambda)t]}$ of type {\rm (4)}.
\end{prop}
\begin{proof}
We 
have that  $f^{-1}(a_1,\ldots,a_{n-1},a_n^{\prime},d) =
(a_1,\ldots,a_{n-1},a_n^{\prime}-\delta,d)$ where   
$$\delta = 
d\ell-\frac{\ell}{\lambda_1}a_1-\cdots-\frac{\ell}{\lambda_{n-1}}a_{n-1}.$$
By assumption 
$$\frac{a_1}{\lambda_1}+\cdots+\frac{a_{n-1}}{\lambda_{n-1}}+\frac{a_n^{\prime}}{\lambda_n+\ell}\geq 
d.$$
Hence 
$$\frac{a_n^{\prime}}{\lambda_n+\ell}\geq 
d-\frac{a_1}{\lambda_1}-\cdots-\frac{a_{n-1}}{\lambda_{n-1}}=
\frac{\delta}{\ell}$$
so $a_n^{\prime}\geq 
(\lambda_n+\ell)\delta/\ell>\delta$ and
$a_n^{\prime}-\delta>0$. 
Furthermore, by $(*)$, 

$$\frac{a_1}{\lambda_1}+\cdots+\frac{a_{n-1}}{\lambda_{n-1}}+\frac{a_n^{\prime}-\delta}{\lambda_n}\geq 
d,$$
hence $(a_1,\ldots,a_{n-1},a_n^{\prime}-\delta,d)$ represents an 
element of 
$\overline{R[I(\blambda)t]}$. \\ Since 
$(a_1,\ldots,a_{n-1},a_n^{\prime},d)$ is minimal, 

$$\frac{a_1}{\lambda_1}+\cdots+\frac{a_{n-1}}{\lambda_{n-1}}+\frac{a_n^{\prime}-1}{\lambda_n+\ell}< 
d.$$
By $(*)$ we may conclude 

$$\frac{a_1}{\lambda_1}+\cdots+\frac{a_{n-1}}{\lambda_{n-1}}+\frac{a_n^{\prime}-1-\delta}{\lambda_n}< 
d.$$
Because $\lambda_n\geq\ell$ (and hence $\lambda_n\geq\lambda_i$ 
for $i<n$), 
we also have that 
\newline 
$(a_1,\ldots,a_{n-1},a_n^{\prime}-\delta,d)$ satisfies condition (B)
for $\overline{R[I(\blambda)t]}$. If condition (C) fails then
write $(a_1,\ldots,a_{n-1},a_n^{\prime}-\delta,d)$ as a sum of at 
least
two minimal generators $(a_{1,i},a_{2,i},...,a_{n,i},d_i)$ of 
$\overline{R[I(\blambda)t]}$. 
If $a_{n,i}>0$, or if $a_{n,i}=0$ and
$$\frac{a_{1,i}}{\lambda_1}+\cdots+\frac{a_{n-1,i}}{\lambda_{n-1}}=d_i$$
then $f(a_{1,i},a_{2,i},...,a_{n,i},d_i)$ is the exponent
vector of an element of $\overline{R[I(\blambda^{\prime})t]}$. If
this holds  for all $i$, then $(a_1,\ldots,a_{n-1},a_n^{\prime},d) = 
f(a_1,\ldots,a_{n-1},a_n^{\prime}-\delta,d)=\sum_i 
f(a_{1,i},a_{2,i},...,a_{n,i},d_i)$, contradicting the minimality of
$(a_1,\ldots,a_{n-1},a_n^{\prime},d)$ as a generator of 
$\overline{R[I(\blambda^{\prime})t]}$.
If, for some $i$, $a_{n,i}=0$ and
$$\frac{a_{1,i}}{\lambda_1}+\cdots+\frac{a_{n-1,i}}{\lambda_{n-1}}>d_i$$
then $f(a_{1,i},a_{2,i},...,a_{n,i},d_i)$ will no longer represent
an element of  $\overline{R[I(\blambda^{\prime})t]}$. But we then have
$$\frac{a_{1,i}}{\lambda_1}+\cdots+\frac{a_{n-1,i}}{\lambda_{n-1}}\geq 
d_i+\frac{1}{\ell}$$
and adding over all $i$ we obtain
\[\frac{a_1}{\lambda_1}+\cdots+
\frac{a_{n-1}}{\lambda_{n-1}}+\frac{a_n^{\prime}-
\delta}{\lambda_n}\geq d+\frac{1}{\ell}.\]
Since $\lambda_n\geq\ell$ it now follows that
\[\frac{a_{1}}{\lambda_1}+\cdots+\frac{a_{n-1}}{\lambda_{n-1}}
+\frac{a_n^{\prime}-\delta-1}{\lambda_n}\geq d, \]
which contradicts the previous observation that
$(a_1,\ldots,a_{n-1},a_n^{\prime}-\delta,d)$ satisfies condition (B)
for $\overline{R[I(\blambda)t]}$. This shows that 
$(a_1,\ldots,a_{n-1},a_n^{\prime}-\delta,d)$ is a minimal generator of
$\overline{R[I(\blambda)t]}$. Finally 
$(a_1,\ldots,a_{n-1},a_n^{\prime}-\delta,d)$ is of type (4) since
$a_n^{\prime}-\delta>0$.
\end{proof}

We can now prove Theorem \ref{congruence}. We remarked after Corollary
\ref{corrange} that there is a bijection between minimal generators
of types (1), (2), (3) for $\overline{R[I(\blambda)t]}$ and 
$\overline{R[I(\blambda^{\prime})t]}$. Proposition \ref{injection} 
gives
an injection of generators of type (4) from 
$\overline{R[I(\blambda)t]}$ to
$\overline{R[I(\blambda^{\prime})t]}$. Since $d$ is preserved under 
these
correspondences the first assertion of Theorem \ref{congruence} 
follows
from Lemma \ref{normalitycriterion}. If $\lambda_n\geq\ell$ then
Proposition
\ref{surjection} gives a bijection of generators of type (4) (hence on
all generators), so the final assertion of Theorem \ref{congruence}
follows again from Lemma \ref{normalitycriterion}. 
\hfill $\square$

\begin{example}\label{revisit} We illustrate the above ideas by 
revisiting
Example \ref{example4}. If $\blambda=(2,3,1)$ then the minimal 
generators
of $\overline{R[I(\blambda)t]}$ are the rows of the following array
(computed with \cite{normaliz}).
\begin{center}$\begin{array}{rrrr}1&0&0&0\\
                    0&1&0&0\\
                    0&0&1&0\\
                    2&0&0&1\\
                    0&3&0&1\\
                    0&0&1&1\\
                    1&2&0&1\end{array}$\end{center}
Here $\ell=6$ and
$f(u_1,u_2,u_3,u_4)=(u_1,u_2,u_3+6u_4-3u_1-2u_2,u_4)$.     There are 
no
generators of type (4). We have $\blambda^{\prime}=(2,3,7)$ and the
minimal generators of $\overline{R[I(\blambda^{\prime})t]}$ are  the 
12
rows of the following array.
\begin{center}$\begin{array}{rrrr}1&0&0&0\\
                    0&1&0&0\\
                    0&0&1&0\\ 
                    2&0&0&1\\
                    0&3&0&1\\
                    0&0&7&1\\
	            1&0&4&1\\
	            0&1&5&1\\
                    1&1&2&1\\
                    0&2&3&1\\
                    1&2&6&2\\	
                    1&2&0&1\end{array}$\end{center}
There are five generators of type (4), namely rows 7 through 11.
The normality of
$I(\blambda^{\prime})$ fails because of the row $(1,2,6,2)$ 
(and Lemma \ref{normalitycriterion}), but 
$I(\blambda)$ is normal. If for example we apply $f^{-1}$ to 
$(1,2,6,2)$ 
we obtain $f^{-1}(1,2,6,2)=(1,2,1,2)$. The vector $(1,2,1,2)$ 
satisfies 
conditions (A) and (B) for $I(2,3,1)$. An expression for
$(1,2,1,2)$ as a sum of minimal generators for $I(2,3,1)$ is 
$(1,2,1,2)=
(1,2,0,1)+(0,0,1,1)$ and $(1,2,6,2)=f(1,2,1,2)=f(1,2,0,1)+f(0,0,1,1)=
(1,2,-1,1)+(0,0,7,1)$, which does not contradict the minimality of
$(1,2,6,2)$ as a generator of $I(2,3,7)$. The argument given in the 
proof
of Theorem \ref{surjection} that $f(1,2,0,1)$ should represent an 
element
of $\overline{R[I(2,3,7)t]}$ does not work since $\lambda_3=1$ is not 
large
enough.
    On the other hand, if we pass to $I(2,3,13)$ we obtain that the
minimal generators of $\overline{R[I(2,3,13)t]}$, by \cite{normaliz},
 are the rows of 
\begin{center}$\begin{array}{rrrr}1&0&0&0\\
                    0&\phantom{1}1&0&\phantom{1}0\\
                    0&0&1&0\\ 
                    2&0&0&1\\
                    0&3&0&1\\
                    0&0&13&1\\
	            1&0&7&1\\
	            0&1&9&1\\
                    1&1&3&1\\
                    0&2&5&1\\
                    1&2&11&2\\	
                    1&2&0&1\end{array}$\end{center}
and indeed $f(1,0,4,1)=(1,0,7,1),\ f(0,1,5,1)=(0,1,9,1),\ 
f(1,1,2,1)=(1,1,3,1),\
f(0,2,3,1)=(0,2,5,1)$ and $f(1,2,6,2)=(1,2,11,2)$, giving the 
bijection
of Proposition \ref{surjection} on generators of type (4) for
$R[I(\blam')t]$ and
$R[I(\blam'')t]$. 
\end{example}
The following example shows that we cannot replace the hypothesis 
$\lambda_n\geq\ell$ by  $\lambda_n\geq\lambda_i$ for all $i$.
However we do not know if the hypothesis $\lambda_n\geq\ell$ is sharp.
\begin{example} $I(2,3,5,6)$ is normal, but $I(2,3,5,36)$ is not.
\end{example}
\begin{remark}
We will conclude with the following remark. Suppose that 
$\lambda_n\geq \lambda_i$ for all $i$.
Then there is a one-to-one correspondence between minimal generators
of $\overline{R[It]}$ with $d=1$ and $(n-1)$-tuples 
$(a_1,\ldots,a_{n-1})$
of non-negative integers $a_i$ such that 
$$\frac{a_1}{\lambda_1}+\cdots+\frac{a_{n-1}}{\lambda_{n-1}}<1$$
under which $(a_1,\ldots,a_{n-1})$ corresponds to 
$(a_1,\ldots,a_{n-1},a_n,1)$, where $a_n$ is the smallest
integer (necessarily positive) such that 
$$\frac{a_1}{\lambda_1}+\cdots+\frac{a_{n-1}}{\lambda_{n-1}}+\frac{a_n}{\lambda_n}\geq
1.$$ 
For $(a_1,\ldots,a_{n-1},a_n,1)$ satisfies condition (B)
because $\lambda_n\geq \lambda_i$ for all $i$, and condition (C) 
because 1 cannot be written as the sum of two positive integers. One can 
attempt to
use this argument to recursively enumerate the minimal generators with
$d=1$.
\end{remark}


\begin{thebibliography}{20}

\bibitem{BG}   Bruns, W.; Gubeladze, J. Rectangular Simplicial
Complexes. In {\em Commutative Algebra, Algebraic Geometry, and 
Computational Methods}; Eisenbud, D., Ed.; Springer: Singapore, 1999; 201-214.


\bibitem{BGT} Bruns, W.; Gubeladze, J.; Trung, N.V. Normal Polytopes,
Triangulations, and Koszul Algebras. J. Reine Angew. Math. {\bf 1997},
{\em 484}, 123-160.

\bibitem{BBV}  Bruns, W.; Vasconceles, V.V.; Villarreal, R.H. Degree bounds in
monomial subrings.  Illinois J. Math. {\bf 1997}  {\em 41}, 341-353.

\bibitem{BH} Bruns, W.; Herzog, J. {\em Cohen-Macaulay Rings}; Cambridge 
University Press: Cambridge, 1993.
 
\bibitem{normaliz} Bruns, W.; Koch, R. NORMALIZ, a Program to Compute
Normalizations of Semigroups. Available by anonymous ftp from
{\bf ftp://ftp.mathematik.Uni-Osnabrueck.DE/pub/osm/kommalg/software/}.


\bibitem{Eis} Eisenbud, D. {\em Commutative Algebra with a View Toward Algebraic
Geometry};  Graduate Texts in Math., Springer-Verlag: New York, 1995; Vol. 150.


\bibitem{F} Faridi, S. Normal Ideals of Graded Rings. Comm. Algebra 
{\bf 2000}, {\em 28}, 1971-1977.


\bibitem{Gil} Gilmer, R. {\em Commutative Semigroup rings}; Chicago Lectures in
Mathematics, The University of Chicago Press: Chicago, 1984.
  



\bibitem{RR} Reid, L.; Roberts, L.G. Monomial Subrings in Arbitrary
Dimension. J. Algebra {\bf 2001}, {\em 236}, 703-730.

\bibitem{RV} L. Reid and M. A. Vitulli, The weak subintegral closure of 
a monomial ideal, Comm. Algebra {\bf 1999}, {\em 27}, 5649-5667.

\bibitem{Rib} Ribenboim, P. Anneaux de Rees Int\'egralement Clos. J. Reine
Angew. Math. {\bf 1960}, {\em 204}, 99-107. 


\bibitem{Zar} Zariski, O.; Samuel, P. {\em Commutative Algebra}; Van Nostrand:
Princeton, 1960; Vol II.



\end{thebibliography}
\end{document}